\documentclass[12pt]{amsart}
\advance\hoffset by -0.5 truecm
\usepackage{amssymb}
\newtheorem{Theorem}{Theorem}[section]

\newtheorem{Definition}[Theorem]{Definition}

\def\V{\mbox{Var}}

\def\R\re
\def\V{\bf V}

\def \re{{\mathbb R}}

\def \0{\lambda_{0}}

\begin{document}
\title[Ricci curvature and Yamabe constants]{Ricci curvature and Yamabe constants}

\author[J. Petean]{Jimmy Petean}\thanks{J. Petean is partially supported 
by grant 46274-E of CONACYT}
 \address{CIMAT  \\
          A.P. 402, 36000 \\
          Guanajuato. Gto. \\
          M\'exico.}
\email{jimmy@cimat.mx}


\date{}


\begin{abstract}We prove that if $(M^n g)$ is a closed Riemannian manifold 
of dimension $n\geq 3$ with volume $V$ and Ricci curvature $Ricci(g) \geq \rho >0$
then the Yamabe constant of the conformal
class $[g]$ satisfies  $Y(M,[g]) \geq n\rho V^{(2/n)}$; 
the equality is achieved if $g$
is an Einstein metric (of Ricci curvature $\rho$). 
This has actually already been proved by S. Ilias \cite{Ilias} in
the context of Sobolev inequalities. 
This implies for instance 
that if $g_1$ is the Fubini-Study metric on ${\bf CP}^2$ and $g$ is any
other Riemannian metric on ${\bf CP}^2$ with $Ricci(g)\geq Ricci(g_1 )$ then
$Vol({\bf CP}^2 ,g) \leq Vol({\bf CP}^2 ,g_1 )$.

\end{abstract}

\maketitle

\section{Introduction}Let $(M^{n},g)$ be a closed Riemannian 
manifold of dimension $n\geq 3$. 
Restricting the total scalar curvature functional to the
conformal class $[g]$ of $g$ we have the Yamabe functional defined
on $L^2_1 (M)$ by

$$Y_g (f) = \frac{a_n \int_M {\| \nabla f \|}^2 dvol_g + \int_M Scal_g \ 
f^2 \ dvol_g }{ {\left( \int_M f^p dvol_g \right) }^{\frac{2}{p}}}.$$

In the expression, and throughtout the article, $a_n = \frac{4(n-1)}{n-2}$,
$p=p_n =\frac{2n}{n-2}$, $dvol_g$ is the volume element of $g$ and 
$Scal_g$ its scalar curvature. 

The Yamabe constant of the conformal class of $g$, $Y(M,[g])$ is the infimum
of this functional. A fundamental result proved in several stages by
Yamabe \cite{Yamabe}, Trudinger \cite{Trudinger}, 
Aubin \cite{Aubin} and R. Schoen \cite{Schoen}
says that there is always a minimizing
function $f_0$ which is smooth and positive. The metric $f_0^{\frac{4}{n-2}} \ 
g$ then has constant scalar curvature and is called a Yamabe metric.

The metric of constant sectional curvature 1, $g_0$, on the sphere 
is a Yamabe metric
and we will denote $Y_n =Y(S^n , g_0 )= n(n-1)V_n^{\frac{2}{n}}$ ($V_n$ is the
volume of $(S^n ,g_0 )$). This value is important in the study of Yamabe
constants since Aubin \cite{Aubin} showed that for any conformal 
class $[g]$ in any closed
$n$-dimensional manifold $M$, $Y(M,[g]) \leq Y_n$
(actually the solution of the Yamabe problem comes from showing that the 
inequality is strict except for the case of $[g_0 ]$). 
It is also easy to check that
$Y(M,[g])$ is positive if and only if there is a metric of positive scalar curvature
on $[g]$ (since the infimum of the Yamabe functional is always realized). One sees
that the study of the Yamabe constant of a conformal class depends strongly on
whether the invariant is positive or non-positive. In the non-positive case
it is particularly useful that the Yamabe constant of the conformal class of a 
metric $g$ is bounded from below by $\inf_M \{Scal_g \} Vol(M,g)^{\frac{2}{n}}$.
This follows by a simple application of H\"{o}lder's inequality to the 
Yamabe functional and it was first pointed out by O. Kobayashi \cite{Kobayashi}. 
This is no longer true in the positive case; by considering Riemannian 
products one can easily build examples of unit volume Riemannian metrics 
with scalar curvature constant and very big ($ \gg Y_n$). 

The aim of this article is to prove that in the positive case there is a
similar lower bound for the Yamabe constant using the infimum
of the Ricci curvature instead of the scalar curvature. Namely, we will prove:

\vspace{1cm}

{\bf Theorem A :} {\it Let $(M^n ,g)$ be a closed Riemannian manifold with
Ricci curvature $Ricci(g) \geq n-1$ and volume $V_0$. Then

$$Y(M,[g])\geq n(n-1)\ V_0^{\frac{2}{n}}
= {\left( \frac{ V_0}{V_n} \right) }^{\frac{2}{n}} \ Y_n .$$ }

\vspace{1cm}

The author was informed by Guofang Wang that the inequality in this Theorem
has already been proved by S. Ilias \cite{Ilias}. Actually the proof given 
in this article goes along the same lines as Ilias' original proof. 
Note that if $g$ is an Einstein metric 
(of constant Ricci curvature n-1) then
it is known to be a Yamabe metric and $Y(M,[g]) = n(n-1)V_0^{\frac{2}{n}}
={\left( \frac{ V_0}{V_n} \right) }^{\frac{2}{n}} \ Y_n $.  The inequality
is therefore optimal. C. B\"{o}hm, M. Wang and W. Ziller \cite{Ziller} have
shown that for $\delta$ close to 1 and $g_0$ the round metric on $S^2$ the
Riemannian metric $\delta g_0 \times g_0$ on
$S^2 \times S^2$ is a Yamabe metric: when 
$\delta \neq 1$ this is probably the simplest case where inequality in
Theorem A is strict. 

\vspace{1cm}

The Yamabe invariant of $M$ was introduced by R. Schoen \cite{Schoen2}
and O. Kobayashi \cite{Kobayashi} as:

$$Y(M) =\sup_{ {[g]}} \ Y(M,[g]),$$

\noindent
the supremum of the Yamabe constants over the space of all conformal classes
of metrics on $M$. Knowledge of the Yamabe invariant and Theorem A produce
some restriction between Ricci curvature and volume of any Riemannian metric
on the given manifold. For instance C. LeBrun \cite{LeBrun} (and M. Gursky and C. LeBrun \cite{Gursky} by more elementary methods) have shown that the Yamabe
invariant of ${\bf CP}^2$ is realized by the conformal class of the
(K\"{a}hler-Eintein) Fubini-Study metric $g_1$. Therefore we obtain:

\vspace{1cm}

{\bf Theorem B :} {\it For any Riemannian metric $g$ on ${\bf CP}^2$ with 
Ricci curvature $Ricci(g) \geq Ricci(g_1 )$ we have $Vol({\bf CP}^2 ,g)
\leq Vol({\bf CP}^2 ,g_1 )$.}

\vspace{1cm}

As another application one recalls that for a Riemannian 4-manifold
$(M,g)$ the space of self-dual 2-forms gives a {\it polarization} of $M$;
namely, a maximal linear subspace of $H^2 (M)$ where the intersection 
form is positive definite \cite{LeBrun2}. 
Now if $g_K$ is a positive K\"{a}hler-Eintein
metric on $M$ and $g$ is any other Riemannian metric on $M$ which defines
the same polarization as $g_K$ then C. LeBrun proved that
\cite[Proposition 2]{LeBrun}  $Y(M,[g]) \leq Y(M,[g_K ])$, and then
we have again that if $Ricci(g) \geq Ricci(g_K )$ then 
$Vol(M,g) \leq Vol(M,g_K )$.

\vspace{1cm}

{\it Acknowledgements:} The author would like to express his gratitude 
to IMPA where this work was carried on with the partial support of
CAPES-Brazil. He would also like to thank Claude
LeBrun for very helpful comments on the original draft of the
manuscript. He would also like to thank Guofang Wang for pointing out
the reference \cite{Illias}.

\section{Spherical rearrangements and isoperimetric inequalities}

In this section we recall a few results we will need for the proof of
Theorem A.

Fix a smooth positive function $f$ on a closed Riemannian manifod $(M,g)$ 
of volume $V_0$. The spherical rearrangement of $f$ is the 
radially symmetric positive
function $f_*$ on $S^n$ such that if we renormalize $S^n$ to have volume 
$V_0$ (and constant sectional curvature) then $\mu (\{ f>t \}) =
\mu (\{ f_* >t \} )$, for all $t\in \re$. 
Here and throughout the article $\mu$ means
the measure corresponding to the volume element of the corresponding
Riemannian metric.

Note that for any positive number $q$ 

$$\int_M f^q =\int_{S^n_{V_0}} f_*^q \ \ .$$

Also recall the coarea formula:

$$ \int_M {\| \nabla f \| }^2 = \int_0^{\infty} \left( \int_{f^{-1} (t) }
\| \nabla f\| d\sigma_t \right) dt,$$

\noindent
and if $t_0$ is a regular value of $f$ then
the function $t\rightarrow \mu (f<t)$ is smooth at $t_0$ and

$$\frac{d}{dt} \mu (f<t) (t_0 )= \int_{f^{-1} (t_0 )} {\|  \nabla f \|}^{-1} 
d\sigma_t $$

\noindent
($d\sigma_t$ means the volume element coming from the induced Riemannian metric).

Let us also recall the following definition introduced by 
B\'{e}rard-Besson-Gallot \cite{Gallot}

\begin{Definition} For any $\beta \in (0,1)$ let $W_{\beta}
=\{ \Omega \subset M: \Omega$ is open with smooth boundary and
$Vol(\Omega )=\beta V_0 \}.$ The isoperimetric function of $(M,g)$
is 

$$h_{(M,g)} (\beta )=h(\beta ) =
\inf \left( \frac{\mu (\partial \Omega )}{V_0} : \Omega \in 
W_{\beta} \right) .$$   
\end{Definition}

B\'{e}rard-Besson-Gallot proved that if the Ricci curvature 
of $(M,g)$, $Ricci(g) \geq n-1$ and $d$ is the diameter 
then $h(\beta ) \geq A(d) \ h_0 (\beta )$; where $h_0$ is the isoperimetric 
function of the round sphere of curvature 1 and

$$A(d)={\left( \frac{\int_0^{\frac{\pi}{2}} \cos^{n-1} (t)\ dt}{
\int_0^{\frac{d}{2}} \cos^{n-1} (t) \ dt} \right) }^{\frac{1}{n}}.$$

Note that $A(d)\geq 1$ by Myers theorem. This is an improvement on
M. Gromov's estimate in \cite[Appendix C]{Gromov}. 
Actually, Gromov's estimate 
(which does not contain the factor $A(d)$) would be enough for the proof of
Theorem A. 
It is well-known that 
$h_0 (\beta ) V_n$ is the area of the $(n-1)$-sphere which bounds a 
geodesic ball of volume $\beta V_n$. 
Note also that if $\lambda$
is a positive constant then the isoperimetric functions of $g$ and
$\lambda g$ are related by $h_{\lambda g} = \frac{1}{\sqrt{\lambda}} h_g$.

\section{Proof of Theorem A}

\begin{proof} Let $f$ be a positive smooth function on $M$
with only non-degenerate (and therefore finite) critical points. 
We will consider
the spherical rearrangement $f_*$ of $f$. We will think of $f_*$ as defined
in the round sphere $S^n_{V_0}$ of volume $V_0$ and therefore 
for any  $t\in \re$, 
$\mu \{ f>t \} =\mu \{ f_* >t \}$.
Note that $S^n_{V_0}$ is obtained by multiplying the round metric of
curvature 1 by ${\left( \frac{V_0}{V_n} \right) }^{\frac{2}{n}}$. One can
put the maximum of $f_*$ in the south pole $q_0$ of $S^n_{V_0}$. Then if
$r$ is the distance in $S^n_{V_0}$ to $q_0$ then $f_*$ is a function of
$r$ and $f_* (r)=t$ if and only if the volume of the geodesic ball of 
radius $r$ in $S^n_{V_0}$ equals $\mu \{ f>t \} $. It follows that if
$t$ is a regular value of $f$ then $f_*$ is differentiable at $r$
and $t$ is a regular value of $f_*$. Note in this case that 
$\| \nabla f_* \|$ is constant along $f_*^{-1} (t)$ since $f_*$ is
radially symmetric. Then we can write

$$ \int_{f_*^{-1} (t) }
\| \nabla f_* \| d\sigma_t  = {\left( \mu (f_*^{-1} (t)) \right)}^2
{\left( \int_{f_*^{-1} (t)} {\| \nabla f_* \|}^{-1} d\sigma_t \right)}^{-1}.$$

We now want to compare the $L^2$-norms of the gradients of $f$ and $f_*$.
By the coarea formula

$$ \int_M {\| \nabla f \| }^2 = \int_0^{\infty} \left( \int_{f^{-1} (t) }
\| \nabla f\| d\sigma_t \right) dt.$$

\noindent
But from H\"{o}lder's inequality (write 
$1={\| \nabla f \|}^{-1/2}{\| \nabla f \|}^{1/2}$)

$$ \int_{f^{-1} (t) }
\| \nabla f\| d\sigma_t  \geq {\left( \mu (f^{-1} (t)) \right)}^2
{\left( \int_{f^{-1} (t)} {\| \nabla f \|}^{-1} d\sigma_t \right)}^{-1}.$$

\noindent
Also note that

$$\int_{f^{-1} (t)} {\| \nabla f \|}^{-1} d\sigma_t = -\frac{d\ \ }{dt}
(\mu \{ f>t \} )$$

$$=-\frac{d\ \ }{dt}
\mu ( \{ f_* >t \} )= \int_{f_*^{-1} (t)} {\| \nabla f_* \|}^{-1} d\sigma_t 
.$$

\noindent
On the other hand $\{ f>t \}$ is a domain in $M$ with volume
$\mu \{ f>t \}$ and boundary $f^{-1} (t)$. By the definition
of the isoperimetric function

$$\mu (f^{-1} (t)) \geq V_0 \ 
h_{(M,g)} \left( \frac{\mu \{ f>t \}}{V_0} \right) .$$

If we let $h_0$ be the isoperimetric function for the sphere then
the estimate of B\'{e}rard-Besson-Gallot says that 

$$h\geq h_0 A(d).$$

The isoperimetric function on the round sphere is realized by round
balls. Therefore

$$\mu (f^{-1} (t)) \geq V_0 \ h_0 \left( \frac{\mu \{ f>t \}}{V_0} \right)
\ A(d) \ 
\ = V_0  
{\left( \frac{V_0}{V_n} \right) }^{\frac{1}{n}}
h_{S^n_{V_0}} \left( \frac{\mu \{ f>t \}}{V_0} \right) \ A(d) $$

$$= {\left( \frac{V_0}{V_n} \right) }^{\frac{1}{n}} \mu ( f_*^{-1} (t) )
\ A(d).$$

And finally,

$$ \int_M {\| \nabla f \| }^2 \geq \int_0^{\infty}  
{\left( \mu (f^{-1} (t)) \right)}^2
{\left( \int_{f^{-1} (t)} {\| \nabla f \|}^{-1} d\sigma_t \right)}^{-1}dt
$$

$$\geq {\left( \frac{V_0}{V_n} \right) }^{\frac{2}{n}} \ (A(d))^2 \ 
\int_0^{\infty}  
{\left( \mu (f_*^{-1} (t)) \right)}^2
{\left( \int_{f_*^{-1} (t)} {\| \nabla f_* \|}^{-1} d\sigma_t \right)}^{-1}dt
$$

$$= {\left( \frac{V_0}{V_n} \right) }^{\frac{2}{n}} \ (A(d))^2 
 \int_0^{\infty} \left( \int_{f_*^{-1} (t) }
\| \nabla f_* \| d\sigma_t  \right)  dt$$

$$= {\left( \frac{V_0}{V_n} \right) }^{\frac{2}{n}} \ (A(d))^2 \ 
\int_{S^n_{V_0}} {\| \nabla f_* \| }^2$$

\noindent
(by the coarea formula).

Therefore

$$Y_g (f)=\frac{a_n \int_M {\| \nabla f \| }^2 + \int_M s_g f^2}{ {\left(\int_M
f^p \right) }^{\frac{2}{p}}} \geq 
\frac{a_n \int_M {\| \nabla f \| }^2 + n(n-1)\int_M f^2}{ {\left(\int_M
f^p \right) }^{\frac{2}{p}}} ,$$

\noindent
since $Ricci_g \geq n-1$. And then from the previous discussion

$$Y_g (f) \geq \frac{a_n V_0^{\frac{2}{n}}
V_n^{\frac{-2}{n}} \ (A(d))^2 \ \int_{S^n_{V_0}} 
{\| \nabla f_* \| }^2 + V_0^{\frac{2}{n}} V_n^{\frac{-2}{n}}
\int_{S^n_{V_0}} Scal_{S^n_{V_0}} f_*^2}{ {\left(\int_{S^n_{V_0}}
f_*^p \right) }^{\frac{2}{p}}}$$

\noindent
And since $A(d)\geq 1$,

$$Y_g (f) \geq {\left( \frac{V_0}{V_n} \right) }^{\frac{2}{n}}
Y_{S^n_{V_0}} (f_* )\geq {\left( \frac{V_0}{V_n} \right) }^{\frac{2}{n}}
Y_n = V_0^{\frac{2}{n}} n(n-1).$$

Since every non-negative function $f\in L_1^2 (M)$ can be approximated 
(in $L_1^2 (M)$)  by a positive Morse function, Theorem A follows by 
taking the infimum for all $f\in L_1^2 (M)$.

\end{proof}

\vspace{0.5cm}

\end{document}